\documentclass[12pt]{amsart}
\usepackage{amsfonts,amssymb,amsmath,amscd}

\usepackage{mathptmx}

\newcommand{\mm}{\mathfrak m}
\newcommand{\nn}{\mathfrak n}
\newcommand{\pp}{\mathfrak p}
\newcommand{\qq}{\mathfrak q}

\newcommand{\N}{\mathbb{N}}

\newcommand{\Z}{\mathbb{Z}}

\newcommand{\Fc}{\mathcal{F}}
\newcommand{\Gc}{\mathcal{G}}

\newcommand{\bsf}{{\boldsymbol f}}

\newcommand{\bsx}{{\boldsymbol x}}

\newcommand{\bsy}{{\boldsymbol y}}

\newcommand{\ges}{\geqslant}
\newcommand{\ov}{\overline}
\DeclareMathOperator{\pnt}{\raise 0.5mm \hbox{\large\bf.}}

\DeclareMathOperator{\ann}{ann}
\DeclareMathOperator{\Coker}{Coker}
\DeclareMathOperator{\depth}{depth}

\DeclareMathOperator{\grade}{grade}
\DeclareMathOperator{\gr}{gr}
\DeclareMathOperator{\Ext}{Ext}
\DeclareMathOperator{\Hom}{Hom}
\DeclareMathOperator{\ild}{inj\, ld}
\DeclareMathOperator{\ilin}{inj\, lin}
\DeclareMathOperator{\indeg}{indeg}
\DeclareMathOperator{\Ker}{Ker}
\DeclareMathOperator{\Image}{Image}

\DeclareMathOperator{\pd}{proj\, dim}
\DeclareMathOperator{\ld}{ld}
\DeclareMathOperator{\gld}{gl\, ld}
\DeclareMathOperator{\lin}{lin}

\DeclareMathOperator{\rank}{rank}
\DeclareMathOperator{\reg}{reg}

\DeclareMathOperator{\Tor}{Tor}

\newcommand{\kos}[2][\bsx]{\operatorname{K}(#1\,;\,#2)}
\newcommand{\pot}[2][k]{#1[\kern-0.28ex[#2]\kern-0.28ex]}
\newcommand{\xra}{\xrightarrow}
\newcommand{\wt}{\tilde}
\newcommand{\shift}{{\scriptstyle\Sigma}}

\newtheorem{thm}{\bf Theorem} [section]
\newtheorem{lem}[thm]{\bf Lemma}
\newtheorem{cor}[thm]{\bf Corollary}
\newtheorem{prop}[thm]{\bf Proposition}

\theoremstyle{definition}
\newtheorem{defn}[thm]{\bf Definition}

\newtheorem{rem}[thm]{\bf Remark}
\newtheorem{ex}[thm]{\bf Example}

\newtheorem{constr}[thm]{\bf Construction}

\theoremstyle{plain}
\newtheorem*{thm*}{Theorem}
\newtheorem{chunk}{}[thm]

\numberwithin{equation}{thm}

\title[Linearity defects]{Linearity defects of modules over commutative rings}

\dedicatory{To Paul Roberts on his sixtieth birthday.}

\author{Srikanth B. Iyengar}
\address{Department of Mathematics, University of Nebraska, Lincoln, NE 68588, USA}
\email{iyengar@math.unl.edu}
\thanks{Part of S.I.'s work was done at the University of Paderborn, on a visit
supported by a Bessel For\-schungs\-preis awarded by the Humboldt-Stiftung.
S.I. was also partly supported by NSF grant, DMS 0602498.}

\author{Tim R\"omer}
\address{Universit\"at Osnabr\"uck, Institut f\"ur Mathematik, 49069 Osnabr\"uck, Germany}
\email{troemer@uos.de}


\begin{document}

\begin{abstract}
This article concerns linear parts of minimal resolutions of finitely generated modules over commutative local, or graded rings.  The focus is on the linearity defect of a module, which marks the point after which the linear part of its minimal resolution is acyclic. The results established track the change in this invariant under some standard operations in commutative algebra. As one of the applications, it is proved that a local ring is Koszul if and only if it admits a Koszul module that is Cohen-Macaulay of minimal degree.  An injective analogue of the linearity defect is introduced and studied. The main results express this new invariant in terms of linearity defects of free resolutions, and relate it to other ring theoretic and homological invariants of the module.
\end{abstract}

\maketitle
%
%
%
\section{Introduction}
\label{intro}
In this article we study linear parts of resolutions of modules over commutative noetherian local, or graded, rings.  Let $R$ be a local ring with maximal ideal $\mm$ and residue field $k$. Any complex $F$ of finitely generated free $R$-modules with $\partial(F)\subseteq \mm F$ has a natural $\mm$-adic filtration;  the associated graded complex with respect to it is denoted $\lin^{R}F$, and is called the \emph{linear part} of $F$.   This construction and invariants derived from it have been investigated by Eisenbud, Fl{\o}ystad, and
Schreyer~\cite{EIFLSC03}, Herzog and Iyengar~\cite{HEIY05}, Okazaki and Yanagawa~\cite{OKYA06}, Yanagawa~\cite{YA06, YA07}, and others.

Let $M$ be a finitely generated $R$-module, or a complex of $R$-modules with $H(M)$ boun\-ded below and degreewise finite, and let $F$ be its minimal free resolution. Herzog and Iyengar~\cite{HEIY05} introduce
the \emph{linearity defect} of $M$ as the number
\[
\ld_R M=\sup\{i \in \Z : H_i(\lin^RF)\neq 0\}.
\]
Following \cite{HEIY05}, a finitely generated $R$-module $M$ is \emph{Koszul} if $\ld_{R}M=0$. Such modules are characterized by the property that their associated graded module $\gr_{\mm}M$ has a linear resolution over the associated graded ring $\gr_{\mm}R$. The ring $R$ is \emph{Koszul} if $k$ is a Koszul module, that is to say, the $k$-algebra $\gr_{\mm}R$ is Koszul, in the classical sense of the word.

We say that $R$ is \emph{absolutely Koszul} if every finitely generated $R$-module has finite linearity defect; equivalently, has a Koszul syzygy module. While absolutely Koszul rings have to be Koszul, the converse does not hold; see the discussion in the introduction of \cite{HEIY05}.
One of the main results of \cite{HEIY05} is that complete intersection local rings and Golod rings are absolutely Koszul. Little else is known about the class
of absolutely Koszul rings.

In Theorem \ref{thm:absolutelykoszul} we prove the following result:

\smallskip

\emph{Let $R\to S$ be a surjective homomorphism of local rings such that
  the projective dimension of the $\gr_{\mm}R$-module $\gr_{\mm}S$ is
  finite. If $S$ is absolutely Koszul, then so is $R$. Moreover, in
  this case, one has an inequality}
\[
\gld R\leq \gld S + \pd_ R S\,.
\]

Here $\gld R$, the \emph{global linearity defect} of $R$, is the supremum of $\ld_RM$ as $M$ ranges over all finitely generated $R$-modules. The proof of the preceding theorem is based on results that track the behavior of linearity defects under some standard operations in commutative algebra: tensor products, quotients by regular sequences, and change of rings. A critical ingredient in the proofs of these latter results is the New Intersection Theorem, in the form of the Amplitude Inequality for complexes. This is the content of Section \ref{sec:bounds}.

A different application of these results concerns the Koszul property of Cohen-Macaulay modules of minimal degree, and is presented in Section \ref{sec:ulrich}.  We say that a Cohen-Macaulay $R$-module $M$ has minimal degree if its degree equals the minimal number of generators of $M$. In Theorem \ref{thm:ulrich:test} we prove that the following statements are equivalent:
\emph{\begin{enumerate}
   \item the ring $R$ is Koszul;
   \item each Cohen-Macaulay $R$-module of minimal degree is Koszul;
   \item there exists a Cohen-Macaulay $R$-module of minimal degree which is Koszul.
\end{enumerate}}

So far our results concern minimal free resolutions of modules (or complexes).  Eisenbud,
Fl{\o}ystad, and Schreyer~\cite{EIFLSC03} considered also minimal injective resolutions over the exterior algebra.  They exploit the fact that over exterior algebras injective modules are free.  Motivated be their results we introduce, in Construction \ref{injective}, a natural filtration on minimal complexes of injective
modules, and the corresponding associated graded complex.  This leads to a notion of the \emph{injective linearity defect} of a module, or a suitable complex, $M$, which we denote $\ild_RM$.

While the definition of the injective linearity defect is straightforward, it is difficult to compute, for minimal injective resolutions are not easily accessible. With this in mind we prove, in Theorem~\ref{thm:duality}, that if the local ring $R$ admits a dualizing complex $D$, suitably normalized, then
\[
\ild_RM=\ld_R \Hom_{R}(M,D)\,.
\]
Thus, one can compute the injective linearity defect using free resolutions, but of the complex $\Hom_{R}(M,D)$. The proof of Theorem~\ref{thm:duality} uses the machinery of local duality theory. One consequence of this result---see Corollary~\ref{cor:injcm}---is an inequality
\[
\ild_{R}M \geq \dim M\,.
\]
This is a little surprising, for the `obvious' lower bound is $\depth M$. As another application of Theorem~\ref{thm:duality}, we prove that when $R$ is Gorenstein and $M$ admits a finite free resolution, say $F$, one has an equality:
\[
\ild_{R}M = \dim R + \sup\{n\mid H_{n}(\lin^{R}\Hom_{R}(F,R))\ne 0\}\,.
\]
We also construct examples that show that the estimates above are optimal.

The results on injective linearity defects are all in Section~\ref{injective linear part of a complex}.

Section \ref{lpdM0} concerns graded rings and modules.
The second author proved in his dissertation \cite{RO01}
that if $R$ is a finitely generated standard graded
Koszul $k$-algebra and $M$ is a finitely generated graded $R$-module,
then $M$ is Koszul if and only if $M$ is componentwise
linear as defined by Herzog and Hibi in \cite{HEHI99}.
That proof has not been published and
we present a compact and simplified version of it here.
In an appendix we collect some technical results related to filtrations needed in the paper.

%
%
%

\section{Bounds on the linearity defect}
\label{sec:bounds}

The starting point of the work in this article is
the construction of the `linear part' of a
complex of modules over a local ring $(R,\mm, k)$, recalled below.

We use the following conventions: Any abelian group $V$ graded by $\Z$ has a lower grading
and an upper grading, and we identify these gradings by setting $V_i=V^{-i}$. We set
\[
\inf V =\inf\{i\in \Z \mid V_i\neq 0\} \quad\text{and}\quad
\sup V =\sup\{i\in \Z \mid V_i\neq 0\}\,.
\]
For any integer $n$, we write $V(n)$ for the graded abelian group with $V(n)_i=V_{n+i}$.

\begin{constr}
\label{projective}
 We say that a complex $F$ of finitely generated free $R$-modules is \emph{minimal} if
$\partial_n(F_n) \subseteq \mm F_{n-1}$ for each $n$.  Let $F$ be such a complex. For each
integer $i$, the graded submodule $\Fc^iF$ of $F$ with
\[
(\Fc^iF)_n=\mm^{i-n} F_n\quad\text{for}\quad n\in \Z\,,
\]
satisfies $\partial(\Fc^iF)\subseteq\Fc^iF$, and hence it is a subcomplex of $F$; as usual, $\mm^j=R$ for $j \leq 0$. Since $\Fc^{i+1}F\subseteq \Fc^iF$ for each $i$, these subcomplexes define a filtration on $F$.  The associated graded complex with respect to it is the \emph{linear part} of $F$, and denoted $\lin^RF$.

Set $A=\gr_\mm R$, the associated graded ring of $R$ with respect to the $\mm$-adic
filtration.  By construction $\lin^RF$ is a complex of graded free $A$-modules with
\[
 \lin^R_{n}F=\gr_\mm(F_n)(-n) \cong  A(-n) \otimes_k F_n/\mm F_n\,,
\]
and the matrices of $\lin^RF$ can be described using linear forms.
\end{constr}

Let $M$ be a complex of $R$-modules whose homology is bounded below and degreewise finite.
Then $M$ has a minimal free resolution: a quasi-isomorphism $F\to M$ where $F$ is a
minimal complex of finitely generated free $R$-modules. Such a complex is unique up to
isomorphism of complexes of $R$-modules and satisfies $F_n=0 $ for $n< \inf H(M)$; for
details see, for instance, \cite[\S1]{RO80}. Herzog and Iyengar \cite{HEIY05} introduced the number
\[
\ld_RM=\sup H(\lin^RF) = \sup\{i \in \Z : H_i(\lin^RF)\neq 0\}
\]
and called it the \emph{linearity defect} of $M$. This number is independent of the choice of $F$,
since minimal resolutions are isomorphic as complexes.

As usual, we identify an $R$-module $M$ with a complex concentrated in degree $0$.  With
this convention, a finitely generated $R$-module $M$ is said to be \emph{Koszul} if
$\ld_RM=0$; the ring $R$ is \emph{Koszul} if $\ld_Rk=0$.

The notion of a Koszul module is motivated by the following considerations.

\begin{rem}
\label{rem:koszul:graded}
  The construction of the linear part of a complex can be carried out also over graded
  rings.  In \cite[Remark 1.10]{HEIY05}, it was observed that a standard graded
  $k$-algebra $R$ is Koszul in the sense of the above definition if and only if $R$
  is a Koszul algebra in the classical sense, that is to say, $k$ has a linear resolution
  over $R$.  Moreover, a local ring $(R,\mm,k)$ is Koszul if and only if $\gr_\mm R$ is
  a Koszul algebra.
\end{rem}

The result below bounds $\ld_{R}M$ in terms of (the linearity defect) of its syzygy modules. In this, its behavior differs from both the depth and the dimension of $M$.

\begin{prop}
\label{prop:ldvssup}
Let $M$ be a complex of $R$-modules with $H(M)$ degreewise finite and bounded below, and
$F$ its minimal free resolution. The following statements hold:

\begin{enumerate}
\item
If $H_n(\lin^RF)= 0$, then $H_n(M)=0$.  In particular, $\ld_RM \geq \sup H(M)$ holds.
\item If $s=\sup H(M)$ is finite, then with $W$ the $R$-module
  $H_s(F_{\ges s})$, one has
\[
 \ld_RM = s + \ld_RW\,.
\]
\end{enumerate}
\end{prop}

\begin{proof}
  Let $\hat R$ denote the $\mm$-adic completion of $R$ and set $\hat M= \hat R \otimes_R
  M$.  Recall that $\hat R$ is also a local ring with maximal ideal $ \mm\hat R$ and that
  the natural homomorphism $R \to \hat R$ is faithfully flat.  Observe that $\hat R
  \otimes_R F$ is a minimal free resolution of $ \hat M$ over $\hat R $ and that one has a
  natural isomorphism $\gr_\mm(F) \cong \gr_{\mm \hat R}( \hat R \otimes_R F)$. Moreover, $\sup H(M)=\sup H(\hat M)$.  One may
  thus replace $R$ and $M$ by $\hat R$ and $\hat M$ respectively and assume that $R$
  is complete.

  (a) One has to prove that the following sequence of $R$-modules is exact:
\[
F_{n+1} \to F_n \to F_{n-1}\,.
\]
For each $n$, the filtration $\{\mm^{i-n} F_n\}_{i\in\Z}$ on $F_n$ is exhaustive
and separated, and $F_n$ is complete with respect to it. The sequence above is
compatible with these filtrations and  the induced associated graded sequence is
exact, by hypothesis.  Now apply Proposition \ref{prop:filtered}.

(b)  Set $G= F_{\ges s}$, and note that $H_{i}(G)=0$ for $i>s$. The complex $\shift^{-s} G$ is thus a minimal free resolution of $W$.  Observe that the natural surjective morphism of complexes $F\to G$ yields a surjective morphism $ \lin^{R} F \to \lin^{R} G$, and that this map is bijective in degrees $n\geq s$. Given the inequality in part (a), this implies the middle equality below:
\[
\ld_R M = \sup H(\lin^RF) = \sup H(\lin^RG) = s + \ld_R W\,.
\]
The other equalities hold by definition.
\end{proof}

The following theorem is one of the main results in this section.

\begin{thm}
 \label{thm:tensors}
 Let $R$ be a local ring, and $M,N$ complexes of $R$-modules with homology degreewise
 finite and bounded below, with minimal free resolutions $F$ and $G$ respectively.
    \begin{enumerate}
    \item When $\pd_RN$ is finite, one has inequalities
\[
\ld_RM + \pd_RN \geq \ld_R(F\otimes_RG) \geq \ld_RM + \inf H(N)\,.
\]
\item When $R$ is regular, then the inequality to the right can be
improved to
\[
\ld_R(F\otimes_RG) \geq \ld_RM + \ld_RN\,.
\]
\end{enumerate}
In particular,  if $\pd_RN$ is finite, then
$\ld_R(F \otimes_R G)<\infty$ if and only if $\,\ld_RM<\infty$.
\end{thm}

The inequality on the right in (a) may fail when $\pd_RN$ is not finite:

\begin{ex}
\label{ex:counterex}
  Let $k$ be a field and $R=\pot{x,y}/(x^2,xy)$. Let $F$ be the complex of $R$-modules
  $0\to R\xra{y}R\to 0$, with the non-zero modules in degrees $0$ and $1$, and $G$ the minimal resolution of the $R$-module $R/Rx$. One has that
\[
\ld_R(F \otimes_R G) = 0 \qquad\text{and} \qquad \ld_R F = 1
\,.
\]
Indeed, $F\otimes_RG\simeq k$, since $y$ is a non-zero-divisor on $R/Rx$. The equality on the left now follows, since the ring $R$ is Koszul. The equality on the right holds by inspection.
\end{ex}

The main ingredient in the proof of Theorem~\ref{thm:tensors}, and also Proposition~\ref{prop:changeofrings} below, is Iversen's Amplitude Inequality~\cite{IV77}, which is an equivalent form of Paul Robert's New Intersection Theorem.  We need versions for unbounded complexes established by Foxby and Iyengar~\cite{FOIY02}, and by Dwyer, Greenlees, and Iyengar~\cite{DWGRIY06}. These are recalled below, in a form convenient for their intended applications.

\begin{rem}
\label{rem:nit}
Let $k$ be a field and $A=\bigoplus_{i\ges 0}A_{i}$ a graded commutative noetherian ring with $A_{0}=k$.
Let $Y$ be a minimal complex of finitely generated graded free $A$-module with $Y_{i}=0$ for $|i|\gg 0$. Here minimality means that $\partial(Y)\subseteq A_{\ges 1}Y$.

For any complex $X$ of graded $A$-modules with $H(X)$ non-zero, degreewise finite, and bounded below,
the following inequalities hold:
\begin{gather}
\label{eq:nit}
\sup H(X) + \sup\{i\mid Y_{i}\ne 0\} \geq \sup H(X\otimes_{A}Y) \geq \sup H(X) +\inf H(Y)\,. \\
\intertext{If $A$ is a polynomial ring, then the inequality on the right can be improved to:}
\label{eq:nitpoly}
\sup H(X\otimes_{A}Y) \geq \sup H(X) + \sup H(Y)\,.
\end{gather}

Indeed, the inequalities in \eqref{eq:nit} are contained in (the graded analogue) of \cite[Theorem 3.1]{FOIY02}, which in turn calls upon \cite[Theorem 5.1]{IV77}; see also \cite[Theorem 5.12]{DWGRIY06}.

Suppose $A$ is a polynomial ring. In proving \eqref{eq:nitpoly}, one may assume $\sup H(X\otimes_{A}Y)$ is finite. It then follows from \eqref{eq:nit} that $\sup H(X)$ is also finite. The right-exactness of tensor products and Nakayama's lemma implies that
\[
\inf H(X\otimes_{A}Y)=\inf H(X) + \inf H(Y)\,.
\]
Thus, the desired inequality follows from \cite[Theorem 5.1]{IV77}.
\end{rem}

The proof of Theorem~\ref{thm:tensors} uses also the following elementary observation.

\begin{lem}
\label{lintensor}
For complexes $F,G$ as in Theorem~\ref{thm:tensors}, and with $A=\gr_\mm R$, there is an
isomorphism of complexes of $A$-modules
\[
\bigl(\lin^R F\bigr) \otimes_{A} \bigl(\lin^R G\bigr)
\cong
\lin^R (F \otimes_R G).
\]
\end{lem}

\begin{proof}
For each $n$ one has natural isomorphisms of $A$-modules
\begin{eqnarray*}
\big(\lin^R F  \otimes_{A} \lin^R G\big)_n
&=&
\bigoplus_{i+j=n} \lin^R_i F \otimes_{A} \lin^R_jG \\
&\cong&
\bigoplus_{i+j=n} \bigl(A(-i) \otimes_k (F_i\otimes_Rk)) \otimes_{A}
\bigl(A(-j) \otimes_k (G_j\otimes_Rk)\bigr)\\
&\cong&
A(-n) \otimes_k \bigl( \bigoplus_{i+j=n}  (F_i\otimes_Rk)\otimes_k
(G_j\otimes_Rk)   \bigr)\\
&\cong&
A(-n) \otimes_k \bigl( (F \otimes_R G)_n \otimes_R k \bigr)\\
&\cong&
\lin^R_n (F \otimes_R G).
\end{eqnarray*}
We leave it to the reader to check compatibility with differentials.
\end{proof}

\begin{proof}[Proof of Theorem~\emph{\ref{thm:tensors}}]
Set $A=\gr_\mm R$; this is a graded commutative noetherian ring with $A_{0}=k$, a field. The complexes of $A$-modules $ \lin^R F$ and $\lin^R G$ are minimal complexes of finitely generated free $A$-modules with $\lin_{i}^{R}F=0=\lin_{i}^{R}G$ for $i\ll 0$. Since $\pd_R N$ is finite, the complex $\lin^R G$ of $A$-modules is finite free and $\lin^{R}_{i}G =0$ for $i>\pd_R N$. We are thus in the context of Remark~\ref{rem:nit}.

(a)  From \eqref{eq:nit} one gets the desired inequalities:
\[
\sup H(\lin^R F) + \pd_R N \geq \sup H(\lin^R (F \otimes_R G) ) \geq \sup H(\lin^R F) +
\inf H(N)\,.
\]

(b) When $R$ is regular, $A$ a polynomial ring, so \eqref{eq:nitpoly} yields an inequality:
\[
\sup H(\lin^R (F \otimes_R G) ) = \sup H(\lin^R F \otimes_A \lin^R G) \geq \sup H(\lin^R
F) + \sup H(\lin^R G)\,.
\]
This is the desired conclusion.
\end{proof}

The next result is in the same spirit as Theorem~\ref{thm:tensors}; the proof is similar.

\begin{prop}
\label{prop:changeofrings}
Let $(R,\mm,k)$ be a local ring and $R\to S$ a surjective homomorphism of rings such that
the projective dimension of the $\gr_{\mm}R$-module $\gr_{\mm}S$ is finite.

Let  $M$ be a complex with homology degreewise finite and bounded below and let
$F$ be its minimal free resolution. Then  one has inequalities
\[
\ld_RM +  \pd_{\gr_{\mm}R}(\gr_{\mm}S) \geq \ld_S(S\otimes_RF) \geq \ld_RM\,.
\]
\end{prop}

\begin{proof}
  Set $\nn = \mm S$; this is the maximal ideal of the local ring $S$.
  Note that $\gr_{\mm}S\cong \gr_{\nn}S$. It is easy to verify that
  one has an isomorphism
\[
\lin^{S}(S\otimes_{R}F) \cong \lin^{R}F \otimes_{\gr_{\mm}R}  (\gr_{\mm}S)
\]
of complexes of modules over $\gr_{\mm}S$. Then \eqref{eq:nit}  applied with $X=\lin^{R}F$ and $Y$ the minimal free resolution of $\gr_{\mm}S$ over $\gr_{\mm}R$ yields the desired result.
\end{proof}

Observe that the hypothesis in the preceding result involves the projective dimension over the associated graded ring. This is not an oversight, but a necessity, as is demonstrated by the following example.

\begin{ex}
Let $k$ be a field, set $R=\pot{x,y,z}/(x^{2},xy+z^{3})$ and $S=R/Rx$, so that
\[
\gr_{\mm}R = k[x,y,z]/(x^{2},xy,z^{3})\quad\text{and}\quad\gr_{\mm}S=k[x,y,z]/(x,z^{3})\,.
\]
It is easy to verify that $\pd_{R}S=1$ whilst $\pd_{\gr_{\mm}R}(\gr_{\mm}S)=\infty$.

The $R$-module $M=R/Ry$ has minimal free resolution $F:=0\to R\xra{y}R\to 0$, so that
\[
\ld_{S}(S\otimes_{R}F) = 0 \quad\text{while}\quad \ld_{R}M = 1\,.
\]
\end{ex}

\begin{defn}
  We say that the ring $R$ is \emph{absolutely Koszul} if $\ld_R
  M<\infty$ for every finitely generated $R$-module $M$. As in \cite{HEIY05}, the
  \emph{global linear defect} of $R$ is the number
\[
\gld R = \sup\{\ld_RM \mid \text{$M$ a finitely generated $R$-module} \}.
\]
\end{defn}

Evidently, when $R$ is absolutely Koszul, it is a Koszul ring, at least when $R$ is graded, but the converse does not hold; see the discussion in the introduction of \cite{HEIY05}. Koszul complete intersection rings and Koszul Golod rings are absolutely Koszul, by \cite[Corollary 5.10]{HEIY05}; the latter also has finite global linearity defect, by \cite[Corollary 6.2]{HEIY05}.

\begin{thm}
\label{thm:absolutelykoszul}
Let $R$ be a local ring and $R\to S$ a surjective homomorphism of rings such that
the projective dimension of the $\gr_{\mm}R$-module $\gr_{\mm}S$ is finite.

If the ring $S$ is absolutely Koszul, then so is the ring $R$. Moreover, one has an inequality
\[
\gld R \leq  \gld S +  \pd_RS  \,.
\]
\end{thm}

\begin{proof}
Let $M$ be a finitely generated $R$-module, with minimal free resolution $F$. Since the projective dimension of $\gr_{\mm}S$ over the ring $\gr_{\mm}R$ is finite, the projective dimension of $S$ over $R$ is finite; see, for example,
\cite[Corollary A3.23]{EI95}. Since $H(S\otimes_{R}F)$ is isomorphic to $\Tor^{R}(S,M)$, one deduces that
\[
s = \sup H(S\otimes_{R}F) \leq \pd_{R}S < \infty\,.
\]
Set $W=H_{s}(S\otimes_{R}F)$. Proposition~\ref{prop:ldvssup}(b) then gives the equality below:
\[
\ld_{R}M  \leq \ld_{S}(S\otimes_{R}F) = \ld_{S} W + s \leq \ld_{S}W + \pd_{R}S\,.
\]
The inequality on the left is by Proposition~\ref{prop:changeofrings}.

When $S$ is absolutely Koszul, the inequalities above yield that $\ld_{R}M$ is finite. Since $M$ was arbitrary, one obtains that $R$ is absolutely Koszul, and moreover that
\[
\gld R\leq \gld S + \pd_{R}S\,.
\]
This completes the proof of the theorem.
\end{proof}

Next we focus on a special case of Theorem~\ref{thm:tensors} where $N$ is a Koszul
complex, for this is the one that is used in the sequel.

\begin{rem}
\label{rem:koszul}
Let $\bsx=x_1,\dots,x_c$ be elements in a commutative ring $R$ and $ \kos R$ the Koszul
complex on $\bsx$; see \cite{BRHE98}.  Given a complex $C$ of $R$-modules, we set
\[
\kos C = \kos R\otimes_R C\,.
\]

Let now $(R,\mm,k)$ be a local ring and $\bsx=x_1,\dots,x_c$ elements in $\mm$. The Koszul
complex $\kos R$ is then a finite free complex of length $c$, hence, for any complex $M$
with homology degreewise finite and bounded below, Theorem~\ref {thm:tensors} yields
inequalities
\[
\ld_R M +c \geq \ld_R \kos M \geq \ld_R M\,.
\]
It should be noted that the Amplitude Inequality, which is the
crucial input in the proof
of Theorem~\ref{thm:tensors}, has an elementary proof when $N$ is the
Koszul complex: one
uses a standard induction argument on $c$ and Nakayama's lemma.
\end{rem}

More precise results are available when $M$ is a module:

\begin{thm}
\label{thm:modx}
Let $(R,\mm,k)$ be a local ring and set $A=\gr_\mm R$.  Let $\bsx=x_1,\dots,x_c$ be
elements in $\mm$, and let $\ov{\bsx}$ be their images in $A_1= \mm/\mm^2$.

The following statements hold for each finitely generated $R$-module $M$.
\begin{enumerate}
    \item If $\bsx$ is regular on $M$ and $M$ is Koszul, then
\[
\ld_R(M/\bsx M)= c-\depth_A(A\ov{\bsx}\,;\,\gr_\mm M)\,.
\]
In particular, $M/\bsx M$ is Koszul if and only if $\,\ov{\bsx}$ is regular on $\gr_\mm
M$.
\item If $\bsx$ is regular on $M$, and $M/\bsx M$ is Koszul, then $M$ is Koszul.
\item
If $\ov{\bsx}$ is regular on $\gr_\mm M$, the $R$-modules
$M$ and $M/\bsx M$ are Koszul simultaneously.
    \end{enumerate}
\end{thm}

\begin{proof}
  When $\ov{\bsx}$ is regular on $\gr_\mm M$, the sequence $\bsx$ is regular on $M$; this
  can be deduced from Proposition~\ref{prop:ldvssup}. Thus, in the rest of the proof we
  may assume that the latter condition holds, and hence that the natural map $\kos M\to
  M/\bsx M$ is a quasi-isomorphism.

  (a) Let $F$ be a minimal free resolution of $M$ over $R$. The quasi-isomorphism $F\to
  M$ then induces a quasi-isomorphism $\kos F\to K(\bsx\,; \,M)$, since $\kos R$ is a
  finite free complex. This gives the first equality below:
\begin{align*}
    \ld_R(M/\bsx M) &=\ld_R \kos F\\
    &= \sup H\big(\lin^R( \kos R\otimes_RF)
    \big) \\
    & = \sup H\big(K(\ov{\bsx}\,;\,A)\otimes_A\lin^RF
    \big) \\
    & = \sup H\big(K(\ov{\bsx}\,;\,A)\otimes_A\gr_\mm M
    \big) \\
    & = c-\depth_A(A\ov{\bsx}\,;\,\gr_\mm M).
\end{align*}
The third one holds by the isomorphism observed in Lemma \ref{lintensor}.
Since $M$ is
Koszul, the map $\lin^RF\to \gr_\mm M$ is a quasi-isomorphism, by \cite[Proposition
1.5]{HEIY05}. It induces a quasi-isomorphism
\[
  K(\ov{\bsx}\,;\,A)\otimes_A\lin^RF\to K(\ov{\bsx}\,;\,A)\otimes_A\gr_ \mm M\,.
\]
This justifies the fourth of the displayed equalities above; the last one holds by
definition.

(b) This follows from Theorem~\ref{thm:tensors}(a) applied with $N=K (\bsx\,;\,R)$.

(c) follows from (a) and (b).
\end{proof}

%

\begin{rem}
\label{rem:ldkoszul}
  The argument for part (a) of the preceding result applies to any complex $M$ with $H(M)$
  degreewise finite and bounded below to yield an equality
\[
\ld_R \kos M = c-\depth_A(A\ov{\bsx}\,;\,\lin^RF)\,.
\]
In particular, with $M=R$ one obtains that
\[
\ld_R \kos R = c-\depth_A(A\ov{\bsx}\,;\,\gr_{\mm}R)\,,
\]
but this can be seen directly. Note that when $\bsx\subseteq \mm^{2}$, one gets $\ld_{R}\kos R = c$.
\end{rem}

\section{Modules of minimal degree}
\label{sec:ulrich}
In this section we apply the results of Section \ref{sec:bounds} to modules of minimal degree, as defined below. We begin by recalling some classical invariants from commutative algebra.

Let $(R,\mm,k)$ be a local ring and $M$ a finitely generated $R$-module.  We write
$\ell_R M$ for the length of an $R$-module $M$, and $\nu_RM$ for its minimal number of
generators; thus one has $\nu_RM=\ell_R(M/\mm M)$.  As is well-known, the following limit exists:
\[
 d!\, \lim_{n\to \infty} \frac{\ell (M/\mm^n M)}{n^d}\quad
   \text{where}\quad {d=\dim M}\,.
\]
This is the \emph{degree} (sometimes referred to as the multiplicity) of $M$, and
denoted $\deg M$.

The following lower bound for the degree is well-known; we sketch an argument for lack of
a suitable reference.

\begin{lem}
  \label{lem:deg}
  If $M$ is a Cohen-Macaulay module over a local ring $R$, then
  $\deg M\ge \nu_RM$.
\end{lem}
\begin{proof}
  This inequality is evident when $\dim M=0$ so suppose $\dim M\ge 1$. Replacing $R$ by
  $R/\ann_RM$ one may assume that $\dim M=\dim R$. A standard argument allows one to
  assume that $k$ is infinite, and then one can find a superficial element $x\in \mm$, not
  contained in any minimal prime ideal of $R$, that is a non-zero-divisor on $M$; see
  \cite[Corollary 8.5.9]{HUSW06}. It then follows from \cite[Proposition 11.1.9]{HUSW06}
  that $\deg M = \deg (M/xM)$. Since $\nu_RM=\nu_R(M/xM)$ holds, an iteration gives the
  desired inequality.
\end{proof}

\begin{defn}
\label{defn:ulrich}
We say that a Cohen-Macaulay module $M$ over a local ring $R$ has \emph{minimal degree} if  $\deg_RM=\nu_RM$ holds.
\end{defn}

Observe that if $Q\to R$ is a surjective homomorphism of local rings, then a Cohen-Macaulay $R$-module $M$ has minimal degree as an $R$-module if and only if it has minimal degree when viewed as a $Q$-module.

When $R$ itself is Cohen-Macaulay, the \emph{maximal Cohen-Macaulay} modules of minimal degree are precisely the Ulrich modules; see the articles of Backelin and Herzog~\cite{BAHE89}, and also that of Brennan, Herzog, Ulrich~\cite{BRHEUL87}, and Ulrich~\cite{UL84}. While it is an open question whether Ulrich modules exist over all Cohen-Macaulay rings, Cohen-Macaulay modules of minimal degree exist over any local ring: $k$ is one such.

We are interested in the linearity of free resolutions of modules of minimal degree. First though we
establish some result, extending those in \cite{HEKU87} for the case when they have maximal dimension and $R$ is Cohen-Macaulay.

\begin{prop}
\label{prop:ulrich:properties}
Let $(R,\mm,k)$ be a local ring, $M$ a Cohen-Macaulay $R$-module of minimal degree, and set $e=\deg_RM$.
\begin{enumerate}
\item When $\dim M=0$, then $M\cong k^e$.  When $\dim M\ge 1$ and $k$ is infinite, there
  exists a superficial $M$-regular sequence $\bsx$ in $\mm\setminus \mm^2$ such that
  $M/\bsx M\cong k^e$.
\item
The $\gr_\mm R$-module $\gr_\mm M$ is Cohen-Macaulay of minimal degree.
\end{enumerate}
\end{prop}

\begin{proof}

  (a) When $\dim M=0$, one has equalities
\[
\ell_RM = \deg_RM= \nu_RM = \ell_R(M/\mm M)\,,
\]
where the second equality holds since $M$ has minimal degree. Thus, $\mm M=0$ and $M\cong k^e$.

Suppose $\dim M\ge 1$ and $k$ is infinite. Arguing as in the proof of Lemma~\ref{lem:deg},
one can construct a superficial $M$-regular sequence $\bsx$ with $\deg_R(M/xM)=\deg_RM$;
one can also ensure that it is in $\mm\setminus \mm^2$, by \cite[Proposition
8.5.7]{HUSW06}. The following equalities then hold:
\[
\deg_R(M/\bsx M) = \deg_RM = \nu_R M = \nu_R(M/\bsx M)\,.
\]
Therefore, $M/\bsx M$ is a zero-dimensional module with the same degree as $M$, and
hence it is isomorphic to $k^e$.

(b) By passing to the $\mm$-adic completion of $R$ if necessary, one can assume that
there exists a regular local ring $(S,\nn,k)$ and a surjective local homomorphism $S \to
R$. Clearly, $M$ has minimal degree also as an $S$-module and $\gr_\nn M\cong \gr_\mm M$ as
$\gr_{\nn}S$-modules.  Replacing $S$ by $R$ one may thus assume that the ring $R$ is
regular.

Choosing an $M$-regular sequence $\bsx$ as in (a) gives the first equality:
\[
\ld_{R}(M/\bsx M) = \ld_R(k^e) = \ld_{R}k =0\,.
\]
The last equality holds because regular local rings are Koszul. Therefore, $\ld_{R}M=0$, that is to say, $M$ is a Koszul module by Theorem~\ref{thm:modx}. Thus, if $F$ is a minimal free resolution of $M$ over $R$, then $\lin^R F$ is a minimal free resolution $\gr_\mm M$ over the ring $A=\gr_{\mm}R$, by \cite[Proposition
1.5]{HEIY05}. This yields an equality $\pd_A(\gr_\mm M)=\pd_{R}M$, and hence the following
(in)equalities hold:
\[
\dim_R M = \depth_R M = \depth_{A} (\gr_\mm M) \leq \dim_{A}(\gr_\mm M) = \dim_R M\,.
\]
The second one is by the Auslander-Buchsbaum Equality.  Hence equality holds in the
middle, that is to say, the $A$-module $\gr_\mm M$ is Cohen-Macaulay.  Since
\[
\deg_A(\gr_\mm M)=\deg_RM\quad\text{and}\quad \nu_A(\gr_\mm M)=\nu_RM
\]
always hold, the $A$-module $\gr_\mm M$ also has minimal degree.
\end{proof}

The gist of the next result is that Cohen-Macaulay modules of minimal degree detect the Koszul property of the ring; see Remark~\ref{rem:herzog} for further comments and  antecedents.

\begin{thm}
\label{thm:ulrich:test}
Let $R$ be a local ring. The following conditions are equivalent:
   \begin{enumerate}
   \item the ring $R$ is Koszul;
   \item each Cohen-Macaulay $R$-module of minimal degree is Koszul;
   \item there exists a Cohen-Macaulay $R$-module of minimal degree which is Koszul.
   \end{enumerate}
\end{thm}

\begin{proof}
  Let $M$ be a Cohen-Macaulay $R$-module of minimal degree; for example $k$, the residue field of $R$. The desired equivalences follow once we prove that  $M$ is a Koszul module if and only if the ring $R$ is Koszul, that is to say, $k$ is a Koszul module.

  We may assume that $k$ is infinite. By Proposition~\ref{prop:ulrich:properties}(a),
  there exists a superficial $M$-regular sequence $\bsx$ in $\mm\setminus \mm^2$ with
  $M/\bsx M\cong k^e$; here $\mm$ is the maximal ideal of $R$. Observe that the image of
  $\bsx$ in $\mm/\mm^2$ is regular on $\gr_\mm M$, since the latter is a Cohen-Macaulay
  module over $\gr_\mm R$, by Proposition~\ref{prop:ulrich:properties}(b).  It is now
  immediate from Theorem~\ref{thm:modx}(c) that $M$ is a Koszul module if and only
  if $k$ is a Koszul module.
\end{proof}

\begin{rem}
\label{rem:herzog}
Let $(R,\mm,k)$ be a local ring and $M$ a finitely generated
$R$-module. Theorem~\ref{thm:ulrich:test} implies the following
statements:
\begin{enumerate}
\item
When $M$ is a Cohen-Macaulay of minimal degree, for \emph{any} surjective homomorphism
$Q\to R$ where $(Q,\qq,k)$ is a Koszul local ring, $M$ is Koszul as
a $Q$-module, since $M$ is also has minimal degree over $Q$.  Thus, the
$\gr_{\qq}Q$-module $\gr_{\qq}M$ has a linear resolution.
\item
If there exists \emph{some} surjective homomorphism $Q\to
R$, where $(Q,\qq,k)$ is Koszul and the $\gr_{\qq}Q$-module $\gr_{\qq}M$
has a linear resolution, then $M$ is Cohen-Macaulay of minimal degree.
\end{enumerate}
In this way, Theorem~\ref{thm:ulrich:test} generalizes
the equivalence (i) $\iff$ (iii) in \cite[Proposition 1.5]{BRHEUL87}.
\end{rem}

%
%
%
%

\section{Injective linear part of a complex}
\label{injective linear part of a complex}

In this section we introduce a notion of an `injective linearity
defect' of a module, and establish results that permit one to compute
it in some cases.

As always,  $(R,\mm,k)$ denotes a local ring.

\begin{constr}
\label{injective}
Let $I$ be a \emph{minimal complex} of injective modules,
that is to say, $I$ is a complex of injective $R$-modules
\[
 \cdots\to I^{n-1} \xra{\partial^{n-1}} I^{n} \xra{\partial^{n}} I^{n+1}\to \dots
\]
with the property that $\Ker(\partial^{n}) \subseteq I^n$ is an injective envelope for each $n\in\Z$.
For each integer $j$ we consider the graded submodule $\Gc^{j}I$ of $I$ with
\[
(\Gc^jI)^n= (0:_{I^n}\mm^{j-n})\,.
\]
The minimality of $I$ implies that $(0:_{I^n}\mm )$, the socle of $I^{n}$, is
contained in $\Ker(\partial^{n})$. It follows, by a straightforward induction on $j$, that the differential $\partial$ of $I$ satisfies:
$$
\partial (\Gc^jI)^n = \partial^n(0:_{I^n}\mm^{j-n}) \subseteq (0:_{I^{n+1}}
 \mm^{j-(n+1)}) = (\Gc^jI)^{n+1}\,.
$$
Therefore $\Gc^{j}I$ is a subcomplex of $I$; note also that $\Gc^{j}I\subseteq \Gc^{j+1}I$.
Hence $\{\Gc^j I\}_{j \in \Z}$ is an \emph{increasing} filtration of the complex $I$.
We call the associated graded complex the \emph{injective linear part} of
$I$, and denote it $\ilin_RI$.
\end{constr}

The injective linear part of $I$ depends only on its $\mm$-torsion subcomplex. This is made precise in the  result below, which is useful for computations. In what follows, given a complex $N$, we write $\Gamma_{\mm}N$ for subcomplex of $\mm$-torsion elements; thus, $(\Gamma_{\mm}N)_{i}=\Gamma_{\mm}(N_{i})$.

\begin{lem}
\label{lem:lingamma}
If $I$ is a minimal complex of injective $R$-modules, then so is the subcomplex $\Gamma_{\mm}I$, and the natural inclusion $\Gamma_{\mm}I\subseteq I$ induces an isomorphism
\[
\ilin_{R}(\Gamma_{\mm}I)\cong \ilin_{R}I
\]
of complexes of $\gr_{\mm}R$-modules.
\end{lem}

\begin{proof}
It follows from the structure theory of injective modules that the subcomplex $\Gamma_{\mm}I$ consists of the injective hulls of $k$ occurring in $I$. It is also easily seen that $\Gamma_{\mm}I$ is a minimal complex.  Thus, the canonical inclusion $\Gamma_\mm I \to I$ induces, for each $j$, morphisms
\[
\Gc^j(\Gamma_\mm I)\to \Gc^j(I)
\]
of complexes of $R$-modules. Since $(0:_{I^{n}}\mm^{j-n})\subseteq \Gamma_{\mm}(I^{n})$, these morphisms are bijective, and hence so is the induced morphism of associated graded complexes; thus, one has an isomorphism
$\ilin_R(\Gamma_\mm I)\cong \ilin_R I$ of complexes of $\gr_{\mm}R$-modules, as desired.
\end{proof}

Each complex $M$ of $R$-modules admits a quasi-isomorphism $M\to I$ where $I$ is a minimal complex of injectives. Such a \emph{minimal injective resolution} is unique up to isomorphism of complexes,
and satisfies $I^{j}=0$ for $j<\inf\{n\mid H^{n}(M)\ne 0\}$; see \cite[\S1]{RO80}.

\begin{defn}
Let $I$ be a minimal injective resolution of a complex $M$. We set
\[
\ild_RM=\sup\{i \in \Z : H^i(\ilin_{R}I)\neq 0\}
\]
and call it the \emph{injective linearity defect} of $M$; this is independent of the choice of $I$.

A module $M$ is \emph{injectively Koszul} if $\ild_{R}M=0$.
\end{defn}

With the hindsight provided by Corollary~\ref{cor:injkoszul},
we remark that $k$ itself is injectively Koszul if and only if it is Koszul, that is to say, $R$ is a Koszul ring.

To each $R$-module $M$, we associated a graded $\gr_{\mm}R$-module denoted $\gr_{\Gc}M$, which in degree $-i$ is the $k$-vector space
\[
(\gr_{\Gc}M)_{-i} = \frac{(0:_{M}\mm^{i+1})}{(0:_{M}\mm^{i})}\,.
\]
Thus, this graded vector space is concentrated in non-positive degrees. Since one has an inclusion $\mm (0:_{M}\mm^{i+1}) \subseteq (0:_{M}\mm^{i})$, there is a natural $\gr_{\mm}R$ action on $\gr_{\Gc}M$, with
\[
(\gr_{\mm}R)_{j} \cdot (\gr_{\Gc}M)_{i}\subseteq (\gr_{\Gc}M)_{i+j}\,.
\]
In other words, $\gr_{\Gc}M$ is a graded module over $\gr_{\mm}R$. Each homomorphism $\varphi\colon M\to N$ of $R$-modules induces a homomorphism of $\gr_{\mm}(R)$-modules $\gr_{\Gc}(\Ker \varphi) \to \Ker(\gr_{\Gc}\varphi)$.

In the result below, $\gr_{\mm}R$ is a graded $R$-module via the surjection $R\to k$, and $\Hom$ denotes the  graded module of homomorphisms.

\begin{lem}
\label{lem:grE}
With $E$ the injective hull of the $R$-module $k$, one has isomorphisms
\[
\gr_{\Gc}E \cong \Hom_{R}(\gr_{\mm}R,E)\cong \Hom_{k}(\gr_{\mm}R,k)
\]
of graded $\gr_{\mm}R$-modules. In particular, $\gr_{\Gc}E$ is the injective hull of $k$ as an $\gr_{\mm}R$-module.
\end{lem}

\begin{proof}
For each $i$, one has an exact sequences of $R$-modules
\[
0\to \frac {\mm^{i}}{\mm^{i+1}} \to  \frac R{\mm^{i+1}}\to \frac R{\mm^i}\to 0\,.
\]
Applying $\Hom_{R}(-,E)$ yields an exact sequence
\[
0\to (0:_{E} \mm^i) \to  (0:_{E}\mm^{i+1}) \to \Hom_{R}(\frac {\mm^{i}}{\mm^{i+1}},E)\to  0\,.
\]
Thus, one has  isomorphisms of $k$-vector spaces
\[
\gr_{\Gc}^{i}E
\cong
\Hom_{R}(\frac {\mm^{i}}{\mm^{i+1}},E)
\cong
\Hom_{k}(\frac {\mm^{i}}{\mm^{i+1}},k)
\]
where the second one holds by adjunction, since $\Hom_{R}(k,E)\cong k$. This yields an isomorphism of graded $k$-vector spaces
\[
\gr_{\Gc}E \cong \Hom_{R}(\gr_{\mm}R,E)\cong \Hom_{k}(\gr_{\mm}R,k)\,.
\]
It is not hard to check that this is compatible with the natural $\gr_{\mm}R$-module structures.
It remains to observe that, by the isomorphism above, $\gr_{\Gc}E$ is the injective hull of $k$ as an $\gr_{\mm}R$-module; see  \cite[Proposition 3.6.16]{BRHE98}.
\end{proof}

 The next result is an analogue of \cite[Proposition 1.5]{HEIY05}.

\begin{prop}
Let $M$ be an $R$-module and $I$ its minimal injective resolution.
\begin{enumerate}
\item
The complex $\ilin_{R}(I)$ consists of direct sums of the injective hull of $k$ over $\gr_{\mm}R$ and is minimal.
\item
The natural map $\gr_{\Gc}M \to H^{0}(\ilin^{R}M)$ is injective;
it is bijective when $M$ is injectively Koszul, and then $\ilin_{R}I$ is a minimal injective resolution of $\gr_{\Gc}M$ over $\gr_{\mm}R$.
\end{enumerate}
\end{prop}

\begin{proof}
(a)  Let $E$ be the injective hull of the $R$-module $k$. For each integer $n$, since $\Gamma_\mm I^n$ is a direct sum of copies of $E$, it follows from Lemma \ref{lem:lingamma} and Lemma~\ref{lem:grE} that
$\ilin_{R}^n I$ is a direct sum of copies of the injective hull of $k$ over $\gr_{\mm}R$.

To verify the minimality of $\ilin_{R}I$, note that one has isomorphisms of complexes
\[
\Hom_{\gr_{\mm}R}(k,\ilin_{R}I) \cong \Hom_{\gr_{\mm}R}(k,\Hom_{R}(\gr_{\mm}R, I))
\cong \Hom_{R}(k,I)\,,
\]
where the first one is a consequence of Lemma~\ref{lem:grE}, and the second one is by adjunction.
The minimality of the complex $I$ implies that the differential on $\Hom_{R}(k,I)$ is zero, and so the same holds for the differential on the complex $\Hom_{\gr_{\mm}R}(k,\ilin_{R}I)$. Hence the complex $\ilin_{R}I$ is minimal, for it consists only of injective hulls of $k$ over $\gr_{\mm}R$.

(b) This follows from (a) and  Proposition~\ref{prop:filtered}(b).
\end{proof}

Observe that $\gr_{\Gc}M$ is non-zero if and only if $\depth_{R}M=0$.
Thus, the preceding result implies that $\depth_{R}M=0$ for any
injectively Koszul module $M$. However, for such a module $\dim M=0$
holds, at least when it is finitely generated. We deduce this from
Corollary~\ref{cor:injcm}, which in turn is obtained from
Theorem~\ref{thm:dualityiso} below. In preparation for stating and proving
the latter result, we recall some properties of dualizing complexes, referring to
Grothendieck~\cite{GR67}, Hartshorne~\cite{HA66} and Roberts~\cite{RO80} for proofs.

\begin{rem}
  \label{rem:dualizing}
 Let $(R,\mm,k)$ be a local ring with a normalized dualizing complex
  $D$. For us, this means that $D$ has the following properties:
\begin{enumerate}
\item $D$ is a minimal complex of injective $R$-modules.
\item $H(D)$ is finitely generated as an $R$-module.
\item $\Ext^0_R(k,D)\cong k$ and $ \Ext^i_R(k,D)= 0 $ for $i \neq 0$.
\end{enumerate}
Up to an isomorphism of complexes, there is only one complex satisfying these properties;
see \cite[Chapter V, \S6]{HA66} and \cite[\S2.2]{RO80} for details. When $R$ is a quotient
of a Gorenstein ring, it has a normalized dualizing complex; see~\cite[??]{HA66}. The converse
result also holds, and was proved by Kawasaki~\cite{KA00}.

Let $M$ be a complex of $R$-modules such that each $H_{i}(M)$ is finitely generated, and set
\[
M^\dag=\Hom_R(M,D)\,.
\]
The following properties of dualizing complexes are used in the sequel.

\begin{chunk}
\label{rem:D0}
One has that $D_{i}$ is a direct sum of injective hulls $E(R/\pp)$, where $\pp$ ranges over all prime ideals with $\dim(R/\pp)=i$. In particular, $D_{i}=0$ for $i\notin [0,\dim R]$.
\end{chunk}

This result is contained in \cite[pp. 58]{RO80}; see also \cite[Chapter V, \S7]{HA66}.

\begin{chunk}
\label{rem:DGor}
Let $J$ be the minimal injective resolution of $R$, viewed as a module over itself.
When the ring $R$ is Gorenstein, $\shift^{d}J$, where $d=\dim R$, is its normalized dualizing complex.
\end{chunk}

See \cite[Chapter V, \S10]{HA66}.

\begin{chunk}
\label{rem:Dquasi}
For any quasi-morphism $M\xra{\simeq} N$ of complexes, the induced map $N^{\dag}\to M^{\dag}$ is also a quasi-isomorphism.
\end{chunk}

This follows from \cite[Chapter II, Lemma 3.1]{HA66}.

\begin{chunk}
\label{rem:Dcoherence}
The $R$-modules $H_{i}(M^{\dag})$ are finitely generated. Moreover, if $H(M)$ is bounded below, respectively, bounded above, then $H(M^{\dag})$ is bounded above, respectively, bounded below.
\end{chunk}

This holds by \cite[Chapter II, Proposition 3.3]{HA66}.

\begin{chunk}
\label{rem:Dduality}
The natural biduality morphism $M\to (M^{\dag})^{\dag}$ is a quasi-isomorphism.
\end{chunk}

When $H(M)$ is bounded, this is \cite[\S2, Theorem 3.5]{RO80}; the general case is contained in \cite[Chapter V, Proposition 2.1]{HA66}.

\begin{chunk}
\label{rem:dimdepth}
When $M$ is a module $\sup H(M^{\dag}) = \dim M$ and $\inf H(M^{\dag}) = \depth M$.
\end{chunk}

This result is a consequence of local duality~\cite[Chapter V, Theorem 6.2,]{HA66} and the Grothendieck Vanishing Theorem~\cite[Theorem 3.5.7]{BRHE98}.
\end{rem}

We require also the following result, for which we could find no suitable reference.

\begin{lem}
\label{lem:fdual}
Assume $H(M)$ is bounded below. Let $F$ be a minimal free resolution of $M$, and $I$ the minimal injective resolution of $M^{\dag}$. With $E$ the injective hull of the $R$-module $k$, one has isomorphisms
\[
\Hom_{R}(F,E)\cong  \Gamma_{\mm}(F^{\dag})\cong \Gamma_{\mm}I
\]
of minimal complexes of injective $R$-modules.
\end{lem}

\begin{proof}
Remark~\ref{rem:D0} implies that $\Gamma_{\mm}D_{0}=E$ and $\Gamma_{\mm}D_{i}=0$ for $i\ne 0$. This gives the isomorphism on the left:
\[
\Hom_{R}(F,E) \cong \Hom_{R}(F,\Gamma_{\mm}D) \cong \Gamma_{\mm}\Hom_{R}(F,D)\,.
\]
The one on the right holds because $D$ is a bounded complex and $F$ is degreewise finite. This justifies the first isomorphism of the Lemma.

It follows from Remark~\ref{rem:Dquasi} that $F^{\dag}$ is an injective resolution of $M^{\dag}$, so one has a homotopy equivalence $I\to F^{\dag}$ of complexes of $R$-modules. This induces a homotopy equivalence $\Gamma_{\mm}I \to \Gamma_{\mm}(F^{\dag})$. Now, both complexes in question are minimal and consist of injectives;
for $\Gamma_{\mm}I$ this is by Lemma~\ref{lem:lingamma}, while for $\Gamma_{\mm}(F^{\dag})$ it holds because it is isomorphic to the complex $\Hom_{R}(F,E)$ which is easily seen to have these properties.
Thus, the morphism $\Gamma_{\mm}I\to \Gamma(F^{\dag})$ must be an isomorphism; see \cite[\S2 Theorem 2.4]{RO80}.
\end{proof}

\begin{thm}
\label{thm:dualityiso}
Let $(R,\mm,k)$ be a local ring with a normalized dualizing complex $D$, and $M$ a complex of $R$-module with $H(M)$ degreewise finite and bounded below. Let $F$ be a minimal free resolution of $M$, and $I$ a minimal injective resolution of $M^{\dag}$.

There exists an isomorphism of complexes of graded $\gr_{\mm}R$-modules
\[
\Hom_{k}(\lin^{R}F,k) \xra{\cong} \ilin_{R}I\,.
\]
\end{thm}

\begin{proof}
Let $E$ be the injective hull of $k$. Lemma~\ref{lem:fdual} gives the first isomorphism below:
\begin{equation}
\label{eq:dtoe}
\ilin_{R}\Hom_{R}(F,E) \xra{\cong} \ilin_{R}(\Gamma_{\mm}I) \xra{\cong} \ilin_{R}I \,.
\end{equation}
The second one is by Lemma~\ref{lem:lingamma}. The filtration $\{\Fc^{i}F\}_{i\ges0}$ of $F$ from Construction~\ref{projective} yields  an exact sequence
\[
0\to \frac{\Fc^{i}F}{\Fc^{i+1}F} \to \frac F{\Fc^{i+1}F} \to \frac F{\Fc^{i}F}\to 0
\]
of complexes of $R$-modules for each $i\ge 0$. This induces the exact sequence in the top row of the
diagram
\[
\begin{CD}
0@>>> \Hom_{R}(\frac F{\Fc^{i}F},E) @>>> \Hom_{R}(\frac F{\Fc^{i+1}F},E)
        @>>> \Hom_{R}(\frac {\Fc^{i}F}{\Fc^{i+1}F},E)@>>> 0 \\
@.   @VV{\cong}V    @VV{\cong}V    @VV{\cong}V      \\
0@>>> \Gc^{i}\Hom_{R}(F,E) @>>> \Gc^{i+1}\Hom_{R}(F,E)
        @>>> \frac{\Gc^{i+1}\Hom_{R}(F,E)}{\Gc^{i}\Hom_{R}(F,E)} @>>> 0
\end{CD}
\]
The isomorphisms on the left and the middle are the natural ones:
\begin{align*}
\Hom_{R}\big(\frac F{\Fc^{i}F},E\big)
    &\cong \bigoplus_{n\in\Z} \Hom_{R}(\frac {F_n}{\mm^{n-i}F_{n}},E) \\
    &\cong \bigoplus_{n\in\Z} \Hom_{R}(\frac R{\mm^{n-i}}\otimes_{R}F_{n},E) \\
    &\cong \bigoplus_{n\in\Z} \Hom_{R}(\frac R{\mm^{n-i}},\Hom_{R}(F_{n},E)) \\
    &=\Gc^{i}\Hom_{R}(F,E).
\end{align*}
The isomorphism on the right, in the ladder of complexes above, thus yields an iso\-mor\-phism of complexes
\[
\Hom_{k}(\lin^{R}F,k) \cong
\bigoplus_{i\in\Z}\Hom_{R}(\frac {\Fc^{i}F}{\Fc^{i+1}F},E) \xra{\cong}
\bigoplus_{i\in\Z}\frac{\Gc^{i+1}\Hom_{R}(F,E)}{\Gc^{i}\Hom_{R}(F,E)} = \ilin^{R}\Hom_{R}(F,E)
\]
The first isomorphism holds because each $\frac {\Fc^{i}F}{\Fc^{i+1}F}$ is a complex of $k$-vector spaces.
Given \eqref{eq:dtoe}, all that is left is to verify that the isomorphism constructed above is compatible with
the $\gr_{\mm}R$-module structures. For this, note that the isomorphism is additive in $F$, so it suffices
to check the compatibility for $F=R$, in which case the map in question is the one from Lemma~\ref{lem:grE},
and $\gr_{\mm}R$-linear.

This completes the proof of the theorem.
\end{proof}

As an first application one obtains the following result, which is reminiscent of the fact that
the Betti numbers (respectively, Bass numbers) of $M$ coincide with the Bass numbers (respectively,
Betti numbers) of $M^{\dag}$; see \cite[\S2, Theorem 3.6]{RO80}. Over Gorenstein rings, it leads to a useful method for computing the injective linearity defect; see Corollary~\ref{cor:gor}.

\begin{thm}
\label{thm:duality}
Let $R$ be a local ring with a normalized dualizing complex $D$. Each complex $M$ of $R$-modules with $H(M)$ degreewise finitely generated has the following properties:
\begin{enumerate}
\item
$\ld_RM=\ild_R(M^\dag)$ when $H(M)$ is bounded below.
\item
$ \ild_RM=\ld_R(M^\dag)$ when $H(M)$ is bounded above.
\end{enumerate}
\end{thm}

\begin{proof}
  (a) Let $F$ a minimal free resolution of $M$ and $I$ a minimal injective resolution of $M^{\dag}$.
  Theorem~\ref{thm:dualityiso} yields the third equality below:
\begin{align*}
\ild_{R}(M^{\dag}) &= \sup\{n\mid H^{n}(\ilin_{R}I)\ne 0\} \\
   &= \sup\{n\mid H^{n}(\Hom_{k}(\lin^{R}F,k))\ne 0\} \\
   &= \sup\{n\mid H_{n}(\lin^{R}F)\ne 0\} \\
   &= \ld_{R}M.
\end{align*}
This gives the desired equality.

(b) When $H(M)$ is bounded above, $H(M^{\dag})$ is bounded below, by Remark~\ref{rem:Dcoherence}, so part (a) yields the second equality below:
\[
\ild_{R}M = \ild_{R}(M^{\dag})^{\dag} = \ld_{R}(M^{\dag})\,.
\]
The first equality holds as $M$ and $(M^{\dag})^{\dag}$ are quasi-isomorphic; see Remark~\ref{rem:Dduality}.
\end{proof}

The other applications of Theorem~\ref{thm:dualityiso} in this section are all via Theorem~\ref{thm:duality}.
The lower bound on $\ild_{R}M$ in the result below holds in full generality; see Corollary~\ref{cor:injcm}.

\begin{cor}
\label{cor:gor}
Let $R$ be a Gorenstein local ring,  $M$ a complex of $R$-modules with $H(M)$ degreewise finitely generated
and $\pd_{R}M$ finite, and $F$ its minimal free resolution.
\begin{enumerate}
\item  One has
\(
\ild_{R}M = \dim R + \sup\{n\mid H_{n}(\lin^{R}\Hom_{R}(F,R))\ne 0\}\,.
\)
\item When $M$ is an $R$-module one then has inequalities
\[
\dim R \geq \ild_{R}M \geq \dim M\,.
 \]
Equality holds on the right when the determinantal ideal
$I_{\nu_{R}M}(\gr_{\mm}(\partial^{F}_{0}))$ in $\gr_{\mm}R$ has grade $0$.
\end{enumerate}
\end{cor}

\begin{proof}
We get the bounds by estimating $\ld_RM^{\dag}$ and applying Theorem~\ref{thm:duality}.

Let $J$ be the minimal injective resolution of $R$, and set $d=\dim R$. Since $R$ is
Gorenstein, $\shift^dJ$, is a normalized dualizing complex; see Remark~\ref{rem:DGor}. One has
then quasi-isomorphisms of complexes:
\[
M^{\dag} = \Hom_{R}(M,\shift^{d}J)\xra{\simeq} \Hom_{R}(F,\shift^{d}J)\xleftarrow{\simeq} \Hom_{R}(F,\shift^{d}R)\cong \shift^{d}\Hom_{R}(F,R)\,.
\]
Since the complex $F$ is finite free and minimal, the same is true of $\shift^{d}\Hom_{R}(F,R)$, so one deduces that the latter is a minimal free resolution of $M^{\dag}$. Therefore one has, by definition, the  first equality below:
\[
\ld_{R}M^{\dag} = \sup H\big(\lin^{R}(\shift^{d}\Hom_{R}(F,R)))
               = d + \sup H\big(\lin^{R} \Hom_{R}(F,R))\,.
\]
This proves (a).

(b) Since $H_{n}(\Hom_{R}(F,R))=\Ext^{-n}_{R}(M,R)$, Proposition~\ref{prop:ldvssup} gives a lower bound:
\[
0 \geq \sup H\big(\lin^{R}\Hom_{R}(F,R)\big) \geq -\grade_{R}M \,.
\]
The upper bound holds because $\Hom_{R}(F,R)_{i}=0$ for $i>0$. Given Theorem~\ref{thm:duality},  the displayed inequalities yield inequalities:
\[
d\geq \ild_{R}M \geq  d -\grade_{R}M = \dim M \,.
\]
The equality holds because $R$ is Cohen-Macaulay. Moreover, equality holds on the right precisely when $H_0(\lin^{R}\Hom_{R}(F,R))\ne 0$ holds.
\end{proof}

The next example demonstrates that Corollary~\ref{cor:gor} is optimal.

\begin{ex}
\label{ex_cor46}
Given non-negative integers $p\ge q\ge r$, there exists a regular local ring $R$ and a finitely generated $R$-module $M$ with
\[
\dim R = p\,,\quad \ild_{R}M=q\,,\quad\text{and}\quad \dim_{R} M = r=\depth_{R}M\,.
\]

Indeed, let $k$ be a field, $\bsx = x_{1},\dots,x_{q}$ and $\bsy = y_{1},\dots,y_{p-q}$
indeterminates over $k$, and set $R=\pot {\bsx,\bsy}$, a power series ring in $\bsx$ and $\bsy$. Choose a regular sequence $\bsf= f_{1},\dots,f_{q-r}$ contained in $(\bsx)^{2}$, and set $M= R/R(\bsf,\bsy)$.
It is clear that $R$ and $M$ have the desired dimension and depth. Now we compute $\ild_{R}M$.

The Koszul complex $\kos[\bsf,\bsy]R$ is a minimal free resolution of $M$ over $R$.
Keeping in mind that $\Hom_{R}(\kos[\bsf,\bsy]R,R)\cong \shift^{r-p}\kos[\bsf,\bsy]R$
one readily obtains
\[
\lin^{R}\Hom_{R}(\kos[\bsf,\bsy]R,R) = \shift^{r-p}\kos[\underline{0},\bsy]A \simeq
    \shift^{r-p}\kos[\underline 0]{A/A\bsy}\,,
\]
where  $A=k[\bsx,\bsy]$, the associated graded ring of $R$, and $\underline 0$ is a sequence consisting of $q-r$ copies of $0$. Therefore Corollary~\ref{cor:gor}(a) yields
\[
\ild_{R}M = p  + r - p + q - r = q\,.
\]
This is the desired result.
\end{ex}

To apply Theorem~\ref{thm:duality} one can often pass to the
completion of the local ring and so ensure the presence of dualizing
complexes. The next result is required for such arguments.

Given a local ring $(R,\mm,k)$ we write $\widehat R$ its $\mm$-adic
completion, and for each complex $M$ of $R$-modules, set $\widehat
M=\widehat R\otimes_{R}M$; this is a complex over $\widehat R$. The
flatness of $\widehat R$ over $R$ entails that when the $R$-module
$H(M)$ is degreewise finite (respectively, bounded below/bounded
above), then the same is true of the $\widehat R$-module $H(\widehat
M)$.  \newcommand{\lch}[2]{{\mathbf{R}\Gamma_{#1}(#2)}}

\begin{prop}
\label{prop:completions}
Let $M$ be a complex of $R$-modules with $H(M)$  degreewise finite.

When $H(M)$ is bounded below $\ld_{\widehat R}(\widehat M) = \ld_{R}M$ holds.

When $H(M)$ is bounded above $\ild_{\widehat R}(\widehat M) = \ild_{R}M$ holds.
\end{prop}

\begin{proof}
  Recall that $\mm \widehat R$ is the maximal ideal of $\widehat R$,
  and that the natural homomorphism
\begin{equation}
\label{eq:assiso}
\gr_{\mm}(R)\to \gr_{\mm \widehat R}(\widehat R)
\end{equation}
of graded $k$-algebras is an isomorphism.

Let $F$ be the minimal free resolution of $M$. Since the $R$-module
$\widehat R$ is flat, the complex $\widehat R\otimes_R F$ is a free
resolution of $\widehat M$ over $\widehat R$; it is evidently also a
minimal one.  Given \eqref{eq:assiso}, it is not hard to verify that
the morphism of complexes $F\to \widehat R\otimes_R F$ induces an
isomorphism
\[
\lin^R(F) \to \lin^{\widehat R}(\widehat R\otimes_RF)\,.
\]
Therefore, the equality $\ld_{R}M=\ld_{\widehat R}(\widehat M)$ holds.

Next we verify the claim about injective linearity defects: Let $M\to I$ and $\widehat M\to J$ be minimal injective resolutions over $R$ and $\widehat R$, respectively.  The morphism $M\to \widehat M$ of complexes of $R$-modules induces a morphism $I\to J$, and hence a morphism
\[
\theta\colon \Gamma_{\mm}I \to \Gamma_{\mm \widehat R}J\,.
\]
This map is a quasi-isomorphism because at the level of homology it is
the homomorphism $H_{\mm}^{\pnt}(M)\to H_{\mm\widehat R}^{\pnt}(\widehat M)$ of local cohomology modules, which is bijective; see \cite[Proposition 3.5.4]{BRHE98}. As the injective hulls of $k$ over $R$ and over $\widehat R$ are isomorphic, one can view both $\Gamma_\mm I$ and $\Gamma_{\mm\widehat R} J$ as complexes of
injectives over $\widehat R$.  These complexes are also minimal, so the quasi-isomorphism $\theta$ is an isomorphism; see \cite[\S2 Theorem 2.4]{RO80}.

The preceding isomorphisms and  Lemma~\ref{lem:lingamma} yield isomorphisms:
\[
\ilin_R I\cong \ilin_R(\Gamma_\mm I)\cong
\ilin_{\widehat R}(\Gamma_\mm J)\cong \ilin_{\widehat R}J\,.
\]
Passing to homology, one gets $\ild_R(M) = \ild_{\widehat R}(\widehat M)$, as desired.
\end{proof}

The following corollary is surprising: given Lemma~\ref{lem:lingamma} it is clear that $\ild_{R}M$ has to be at least $\depth M$; it is a priori not clear why it should be greater than $\dim M$.

\begin{cor}
\label{cor:injcm}
Let $R$ be a local ring and $M$ a finitely generated $R$-module.
The inequality $\ild_{R}M\geq \dim M$ then holds. Hence, if $M$ is injectively Koszul, then $\dim M=0$.
\end{cor}

\begin{proof}
  Thanks to Proposition~\ref{prop:completions}, one may pass to the
  completion of $R$ and assume that it has a dualizing complex.
  Theorem~\ref{thm:duality} then yields the first equality below:
\[
\ild_{R}M = \ld_{R}(M^{\dag})\geq \sup\{i\mid H_{i}(M^{\dag})\ne 0\} =\dim M\,;
\]
the inequality is due to Proposition~\ref{prop:ldvssup}; for the last equality, see \ref{rem:dimdepth}.
\end{proof}

With regards to the preceding result, note that $k$ is
zero-dimensional but $\ild_R(k)=0$ holds if and only if the ring $R$
is Koszul; this is by Corollary~\ref{cor:injkoszul} below.

\begin{cor}
\label{cor:injkoszul}
Let $(R,\mm,k)$ be a local ring. The $R$-module $k$ is injectively Koszul if and only if
the ring $R$ is Koszul.
\end{cor}

\begin{proof}
  Since $\widehat k=k$, one can apply
  Proposition~\ref{prop:completions} to pass to the completion of $R$,
  and thus assume that it has a dualizing complex. Since $k^{\dag}=k$, Theorem~\ref{thm:duality} yields that $\ild_{R}k=0$ if and only if $\ld_{R}k=0$, that is to say, the ring $R$ is Koszul.
\end{proof}

\section{Componentwise linear modules}
\label{lpdM0}

Let $k$ be a field and $R$ a \emph{standard graded} $k$-algebra, that is to say,
$R=\bigoplus_{i\in \N} R_i$ is a graded ring with $R_{0}=k$, $\rank_kR_1$ is finite,
and $R=k[R_1]$. In particular, the ring $R$ is noetherian and
$\mm=\bigoplus_{i\ges 1} R_i$ is
the unique graded maximal ideal. Each finitely generated graded $R$-module $M$ admits a
minimal graded free resolution $F$, and its linear part, $\lin^RF$, is defined as in the
local case; see \ref{projective}.  This gives rise to the invariant $\ld_RM$ and a notion
of a Koszul module. As noted in Remark~\ref{rem:koszul:graded}, the ring $R$ is Koszul
precisely when it is Koszul in the classical sense of the word.

In this section, we present a characterization of Koszul modules over Koszul algebras,
which was first established in the second author's thesis~\cite{RO01}. The argument
presented here is a streamlining of the original one.

\begin{rem}
  Observe that since $R$ is standard graded $\gr_\mm(R(-n))$ is naturally isomorphic to
  $R$.  To be more precise one should view $R$ as a bigraded $k$-algebra with components
\[
R_{p,q}=\begin{cases}
R_p & \text{for $p=q$}\,,\\
0 & \text{for $p\neq q$}\,.
\end{cases}
\]
Now let $M$ be a finitely generated graded $R$-module and $F$ its minimal graded free resolution.
For each integer $n\ge 0$, there is an isomorphism
\[
F_n=\bigoplus_{i\in \Z} R(-i)^{\beta^R_{n,i}(M)}\,,\quad
\text{where}\quad \beta^R_{n,i}(M)=\dim_k \Tor_n^R(k,M)_i\,.
\]
The $\beta^{R}_{n,i}$ are the \emph{graded Betti numbers} of $M$. It is then clear that
\[
\lin^R_nF \cong \bigoplus_{i\in \Z}R(-n,-i)^{\beta_{n,i}^R(M)}\,.
\]
The least degree of a generator of $M$ is called \emph{initial degree} of $M$, and denoted
$\indeg M$. Note that $\indeg M =\min\{t \in\Z\colon M_t \neq 0 \}$.
\end{rem}

\begin{defn}
\label{defn:cmr}
  The \emph{Castelnuovo-Mumford regularity} of $M$ is the number
\[
  \reg_R M = \sup\{r \in \Z \mid \beta^{R}_{n,n+r}(M) \neq 0 \text{ for some } n \in \N \}\,.
\]
Note that $\reg_{R}M\geq \indeg M$, with equality if and only if $M$ has an
\emph{$i$-linear resolution}:
\[
  \beta_{n,r}^R(M)= 0\quad\text{for}\quad r \neq i+\indeg M\,;
\]
equivalently, if the differentials in $F$ can be represented by matrices of linear
forms. The complexes $\lin^RF$ and $F$ are then isomorphic, so $\ld_R M=0$; that is to
say, $M$ is Koszul.
\end{defn}

\begin{defn}
For  each $i \in \Z$ let $M_{\langle i \rangle}$ be the submodule of $M$ generated by $M_i$.
The module $M$ is \emph{componentwise linear} if $M_{\langle i \rangle}$ has an $i$-linear
resolution for each $i$.
\end{defn}

\begin{lem}
   \label{linearwithm}
   Let $R$ be a Koszul algebra and $M$ a finitely generated graded $R$-module.  If $M$ has
   an $i$-linear resolution, then $\mm M$ has an $(i+1)$-linear resolution.
\end{lem}

\begin{proof}
  Since $M$ has an $i$-linear resolution, it is generated in degree $i$. Thus $M/\mm M
  \cong \bigoplus k(-i)$ has an $i$-linear resolution because $R$ is a Koszul algebra.  It
  follows from the exact sequence $0 \to \mm M \to M \to M/\mm M \to 0$ that
\[
i+1 =\indeg(\mm M) \leq \reg_R(\mm M) \leq \max \{i,i+1\}\,.
\]
Thus $i+1 =\reg_R(\mm M)$ and  $\mm M$ has an $(i+1)$-linear resolution.
\end{proof}

\begin{lem}
\label{linear_lpd0_compi}
Let $R$ be a Koszul algebra and $M$ a finitely generated graded $R$-module.
The following statements are equivalent:
   \begin{enumerate}
\item $M$ is componentwise linear;
\item $M/M_{\langle \indeg M \rangle}$ is componentwise linear and $M_{\langle \indeg
    M\rangle}$ has a linear resolution.
   \end{enumerate}
\end{lem}

\begin{proof}
  We may assume $\indeg M=0$. Evidently $M_{\langle 0 \rangle\langle i\rangle} =
  \mm^{i}M_{\langle 0\rangle }$ holds, so the sequence
\[
0 \to M_{\langle 0 \rangle\langle i \rangle} \to M_{\langle i
\rangle} \to (M/M_{\langle 0
   \rangle})_{\langle i \rangle} \to 0
\]
is exact. Moreover, when $M_{\langle 0 \rangle}$ has a $0$-linear resolution, $M_{\langle
  0 \rangle\langle i \rangle}$ has an $i$-linear resolution, by
Lemma~\ref{linearwithm}. The equivalence of (a) and (b) now follows from the sequence
above.
\end{proof}

For the next result we recall that over Koszul algebras the regularity of each finitely
generated module is finite; see \cite{AVEI92} and \cite{AVPE01}.

\begin{thm}
 \label{graded:mainresult}
   Let $R$ be a Koszul algebra and $M$ a finitely generated graded $R$-module.
The module $M$ is Koszul if and only if it is componentwise linear.
\end{thm}

\begin{proof}
  Let $F$ be a minimal graded free resolution of $M$ over $R$. Set $d=\indeg M$ and
  consider a graded submodule $\wt{F}$ of $F$ with
\[
\wt{F}_n=R(-n)^ {\beta^R_{n,n+d}(M)}  \quad\text{for}\quad n\geq 0\,.
\]
Observe that, for degree reasons, $\partial (\wt{F})\subseteq \wt{F}$, where
$\partial$ is the differential on $F$, so $\wt{F}$ is a subcomplex of $F$.
Set $\wt{M}=  H_0(\wt{F})$ and observe that
\begin{equation}
\label{eq:lin0}
\wt M = H_0(\wt{F})_{\langle d\rangle} =
H_0(F)_{\langle d\rangle} \cong M_{\langle d\rangle}\,.
\end{equation}
One has an exact sequence of complexes
\begin{equation}
\label{eq:linses}
0 \to \wt{F} \to F \to F/\wt{F} \to 0
\end{equation}
which, by construction, is split as a sequence of graded-modules. Again, degree
considerations reveal that this induces a decomposition of complexes of $R$-modules:
\begin{equation}
\label{eq:linsplit}
\lin^R F = \lin^R(\wt{F}) \oplus \lin^R (F/\wt{F})\,.
\end{equation}

We induce on $\reg_R M -d$ to prove the desired equivalence.  If $\reg_R M =d$, then $M$
has a linear resolution, and hence it is Koszul, as noted in Definition~\ref{defn:cmr},
and componentwise linear, by Lemma~\ref{linear_lpd0_compi}.  Assume $ \reg_R M-d\geq 1$.

When $M$ is Koszul, so that $H_i(\lin^RF)=0$ for $i\geq 1$, one obtains from \eqref{eq:linsplit}
that
\[
H_i(\lin^R(\wt F))=0=H_i(\lin^R(F/\wt F))\quad\text{for}\quad i\geq 1\,.
\]
Proposition~\ref {prop:ldvssup} then yields $H_i(\wt{F})=0=H_i(F/\wt{F})$ for $\geq 1$.
It then follows from \eqref{eq:lin0} and the homology exact sequence arising from
\eqref{eq:linses} that $\wt{F}$ is the minimal free resolution of $\wt{M}$ and $F/\wt{F}$
is the minimal free resolution of $M/\wt{M}$.  The displayed equalities then imply that
$\wt M$ has a linear resolution and $M/\wt M$ is Koszul.  Observing that $\reg_R M
-d>\reg_R (M/\wt M) -\indeg(M/\wt M)$ the induction hypothesis yields that $M/\wt M$ is
componentwise linear, so $M$ is componentwise linear, by Lemma~\ref{linear_lpd0_compi}.

Assume now that $M$ is componentwise linear; then so are $\wt M$ and $M/\wt M$, by
Lemma~\ref{linear_lpd0_compi}.  Because $\wt M$ has a $d$-linear resolution
one obtains the last equality below:
\[
\rank_R \wt{F}_n = \beta_n^R(M)_{n+d}= \beta_n^R (\wt M)_{n+d}=\beta_n^R (\wt M)\,.
\]
The second equality holds because $\wt M=M_{\langle d\rangle}$. An induction on $n$ then
shows that $\wt{F}$ is the minimal free resolution of $\wt{M}$. Hence \eqref{eq:linses}
implies $F/\wt{F}$ is the minimal free resolution of $M /\wt M$.  The induction hypothesis
yields $H_i(\lin^R \wt{F})=0=H_i(\lin^R (F/\wt{F}))$ for $i\geq1$, so $H_i(\lin^RF)=0$ for
$i\ge 1$, by \eqref{eq:linsplit}. Thus, $M$ is Koszul.
\end{proof}

%
%
%

\appendix

\section{Filtrations}
In this paper we need some facts about filtrations.  For the convenience of the reader we
state these results separately in this appendix and present their proofs.

Let $R$ be a ring.  A \emph{filtered} module $U$ is an $R$-module with filtration
$\{U^n\}_{n \in \Z}$ such that $U^{n+1} \subseteq U^{n}$ for $n \in \Z$.  The filtration
is \emph{separated} if $\bigcap_{n\in \Z} U_n=0$ and it is \emph{exhaustive} if
$\bigcup_{n \in \Z} U_n=U$.  The module $U$ is \emph{complete} with respect to the
filtration if the natural map $U \to\varprojlim_n U/U_n$ is an isomorphism.  The
associated graded module of filtered module $U$ is the graded module $\gr\, U$ with
degree $n$-component $U^n/U^{n+1}$.

The proof of the following lemma can be found in
\cite[Chapter III]{BO89}.
\begin{lem}
\label{lem:filtered}
Let $U$ be an $R$-submodule of a filtered $R$-module $V$. Then
\begin{enumerate}
\item
$U$ is a filtered $R$-module with $U^n=U \cap V^n$.
\item
$V/U$ is a filtered $R$-module with $(V/U)^n=V^n/(U \cap V^n)$.
\item
Considering $U$ and $V/U$ as filtered $R$-modules
induced by the filtrations of (a) and (b) respectively
the associated graded sequence below is exact:
\[
0 \to \gr\, U \to \gr\, V \to \gr(V/U) \to 0\,.
\]
\end{enumerate}
\end{lem}

However, observe that $\gr(\cdot)$ is usually not an exact functor.
\begin{ex}
Let $R=k[x]$ be a polynomial ring over a field $k$, and set $\mm=(x)$.
Applying $\gr_\mm(\cdot)$ to the exact sequence
\[
0 \to k[x] \overset{x^2}{\to} k[x] \to k[x]/(x^2) \to 0
\]
leads to the sequence $0 \to k[x] \overset{0}{\to} k[x] \to k[x]/(x^2) \to 0$,
which is not exact. The problem here is that the filtration of $k[x]$ is not compatible
with the filtration of the image $(x^2)$ of the multiplication map by ``$x^2$''
as a submodule of $k[x]$.
\end{ex}

A homomorphism of filtered modules is an $R$-module homomorphism $ \varphi \colon U \to V$
such that $\varphi(U^n) \subseteq V^n$. Such a map induces a homomorphism $\gr\,\varphi\colon
\gr\, U\to \gr\, V$. It follows from Lemma \ref{lem:filtered}
that $\Ker \varphi $ is a filtered $R$-module with $(\Ker \varphi)^n
= \Ker \varphi \cap U^n$ and $ \Coker \varphi$ is a filtered $R$-module with $(\Coker
\varphi)^n = \bigl(V^n + \varphi(U) \bigr)/\varphi(U)$.

\begin{prop}
\label{prop:filtered}
Let $U \overset{\varphi}{\to} V \overset{\psi}{\to} W$ be a sequence of filtered
$R$-modules. If the associated graded sequence is exact, the following statements hold.
\begin{enumerate}
\item
The canonical homomorphism $\Coker(\gr\,\psi) \to \gr(\Coker \psi )$ is bijective.
\item
The canonical homomorphism $\gr\, \Ker(\varphi)\to \Ker(\gr\,\varphi)$ is injective;
it is bijective when the sequence $U \overset{\varphi}{\to} V \overset{\psi}{\to} W$ is also exact.
\item
When  $U$ is complete and the filtration on $V$ is exhaustive and separated,
the sequence $U \overset{\varphi}{\to} V \overset{\psi}{\to} W$ is exact.
\end{enumerate}
\end{prop}

\begin{proof}
(a) This was proved in \cite[Lemma 1.16]{HEIY05}.

(b) Since one has the following equalities:
\[
\gr( \Ker \varphi)^n=\frac{\Ker \varphi \cap U^n}{\Ker \varphi \cap U^{n+1}}  \quad\text{and}\quad
\Ker(\gr\,\psi)^n=\{ \overline{u}\in U^{n}/U^{n+1} : \psi(u) \in V^{n+1}\}
\]
one deduces that the canonical homomorphism $\gr(\Ker\varphi) \to \Ker(\gr\,\psi)$ is injective.
Applying Lemma \ref{lem:filtered} (c) to the exact sequence
$$
0 \to \Ker \varphi \to U \to (U/\Ker \varphi) \to 0
$$
yields an exact sequence
$$
0 \to \gr(\Ker \varphi) \to U \to \gr(U/\Ker \varphi) \to 0.
$$
Assume now that
$
U \overset{\varphi}{\to} V \overset{\psi}{\to} W$
is exact.
Then $U/\Ker \varphi\cong \Image \varphi = \Ker \psi$ as $R$-modules.
Moreover, this isomorphism is compatible with the induced filtrations on these
modules and we obtain an isomorphism
$$
\gr(U/\Ker \varphi)\cong \gr(\Ker \psi).
$$
The map $\gr\,\varphi$ factors as
$$
\gr\, U \to \gr(U/\Ker \varphi)\cong \gr(\Ker \psi) \hookrightarrow \Ker(\gr\,\psi).
$$
It follows from the assumption that  $\gr\, U \to \Ker(\gr\,\psi)$ is surjective.
Hence $\gr(\Ker \psi) \cong \Ker(\gr\,\psi)$ as desired.

(c)
We have to show that the homomorphism $U \to \Ker \psi$ is
surjective. Applying (b) to $\Ker\psi$ yields that
$\gr(\Ker\psi)$ is a submodule of $\Ker(\gr\,\psi)$. The map $\gr\,\varphi$ factors as
\[
\gr\, U \to \gr(\Ker\psi)\subseteq \Ker(\gr\,\psi)\,.
\]
The hypothesis is that $\gr\, U \to \Ker (\gr\,\psi)$ is surjective, so $\gr\, U \to \gr(
\Ker \psi)$ is surjective.  Since $V$ is exhaustive and separated, the same is true for
$\Ker \psi$, with induced filtration. Now it remains to apply \cite[Chapter III, \S2.8,
Corollary 2]{BO89}.
\end{proof}

\end{document}